\documentclass[conference]{IEEEtran}
\IEEEoverridecommandlockouts

\usepackage{cite}
\usepackage{hyperref}
\usepackage{amsmath,amssymb,amsfonts}
\usepackage[ruled,vlined]{algorithm2e}
\usepackage{graphicx}
\usepackage{subcaption}
\usepackage{textcomp}
\usepackage{xcolor}
\usepackage[capitalize]{cleveref}
\usepackage{amsthm}
\usepackage{booktabs}
\usepackage{multirow}
\usepackage{enumitem}
\newtheorem{defn}{Definition}[section]
\newtheorem{thrm}{Theorem}[section]     
\newtheorem{lem}[thrm]{Lemma}
\newtheorem{prob}{Problem}[section]
\usepackage[braket]{mathtools}  
\newcommand{\R}{\mathbb{R}}
\newcommand{\set}[1]{\{\,#1\,\}}
\newcommand{\bigO}{\mathcal{O}}
\newcommand{\defeq}{\coloneqq}
\DeclarePairedDelimiter{\pqty}{(}{)}
\DeclarePairedDelimiter{\Set}{\{}{\}}
\DeclarePairedDelimiter{\abs}{\lvert}{\rvert}
\newenvironment{pf}{\begin{proof}}{\end{proof}}
\crefname{thrm}{theorem}{theorems}
\Crefname{thrm}{Theorem}{Theorems}   
\crefname{lem}{lemma}{lemmas}
\Crefname{lem}{Lemma}{Lemmas}
\crefname{prob}{problem}{problems}
\Crefname{prob}{Problem}{Problems}
\crefname{defn}{definition}{definitions}
\Crefname{defn}{Definition}{Definitions}
\DeclareMathOperator*{\argmax}{arg\,max}
\SetKwFunction{MSSP}{\textsc{MSSP}}
\SetKwFunction{SeedTree}{\textsc{SeedTree}}
\SetKwFunction{ComputePriority}{\textsc{ComputePriority}}
\SetKwFunction{DetachTarget}{\textsc{DetachTarget}}
\SetKwFunction{ReattachTarget}{\textsc{ReattachTarget}}

\begin{document}

\title{Optimizing a Model-Agnostic Measure of Graph Counterdeceptiveness via Reattachment\\
\thanks{This work was supported by the Office of Naval Research through Research Grant no. N00014--23--1--2651.}
}

\author{%
  \IEEEauthorblockN{Anakin Dey\textsuperscript{1},
                    Sam Ruggerio\textsuperscript{2},
                    Manav Vora\textsuperscript{2},
                    Melkior Ornik\textsuperscript{2}}%
  \\[1ex]
  \IEEEauthorblockA{\textsuperscript{1}The Ohio State University, Columbus, OH 43210, USA\\
                    dey.92@buckeyemail.osu.edu}\\[.5ex]
  \IEEEauthorblockA{\textsuperscript{2}University of Illinois Urbana-Champaign, Urbana, IL 61801, USA\\
                    \{samuelr6, mkvora2, mornik\}@illinois.edu}
}

\maketitle

\begin{abstract}
Recognition of an adversary's objective is a core problem in physical security and cyber defense.
    Prior work on target recognition focuses on developing optimal inference strategies given the adversary's operating environment. 
    However, the success of such strategies significantly depends on features of the environment. 
    We consider the problem of optimal counterdeceptive environment design: construction of an environment which promotes early recognition of an adversary's objective, given operational constraints. 
    Viewed as a bounded-length graph-design problem, we introduce a metric
for counterdeception and a novel heuristic that maximizes it based on iterative reattachment of trees.
    We benchmark the performance of this algorithm on synthetic networks as well as a graph inspired by a real-world high-security environment, verifying that the proposed algorithm is computationally feasible and yields meaningful network designs.
\end{abstract}


\section{Introduction}\label{sec:intro}

In many security–critical domains---military bases, industrial plants,
city-scale transportation networks---an \emph{observer} (defender) must infer the true
goal of a potentially adversarial \emph{agent} (attacker) from the
agent’s motion.  The attacker, in turn, plans its path to hide that goal
for as long as possible, a setting studied under \emph{deceptive path
optimisation} and goal recognition
\cite{dpp_dynamic,deceptive_robot_motion,goal_recognition_in_path,prob_plan_recognition}.
Almost all prior work treats the layout of the facility as immutable; at most, the defender is allowed to disable or label a handful of individual corridors, roads, or communication links \cite{game-theory-goal-recognition,active_goal_rec_design}.
On a fixed map, much of the goal-recognition work performs probabilistic inference under a shortest-path assumption for the adversary
\cite{goal_recognition_design_survey,game-theory-goal-recognition}. Some efforts relax that assumption by allowing the agent to take boundedly sub-optimal or stochastic routes, but still keep the environment fixed
\cite{goal_rec_design_non_opt,goal_rec_design_non_obs}.  Road-network design
studies optimize flow or coverage under budgets
\cite{urban_road_network}, but they do not consider deception; conversely, optimal target placement has been explored without an
explicit road network \cite{dispersion_models,measure_target_predict}.  In many real scenarios the defender \emph{owns} the infrastructure:
roads on a military air-base, hallways in a secure facility, even
routing rules in cyberspace.  The natural question becomes: \emph{How should we design a network so that an attacker’s true
destination becomes obvious as early as possible, subject to a global
``construction'' budget?} 

To our knowledge, full network \emph{design} against a worst-case, goal-hiding
adversary has not been addressed.

\begin{figure}[h]
  \centering
  \includegraphics[width=.2\textwidth]{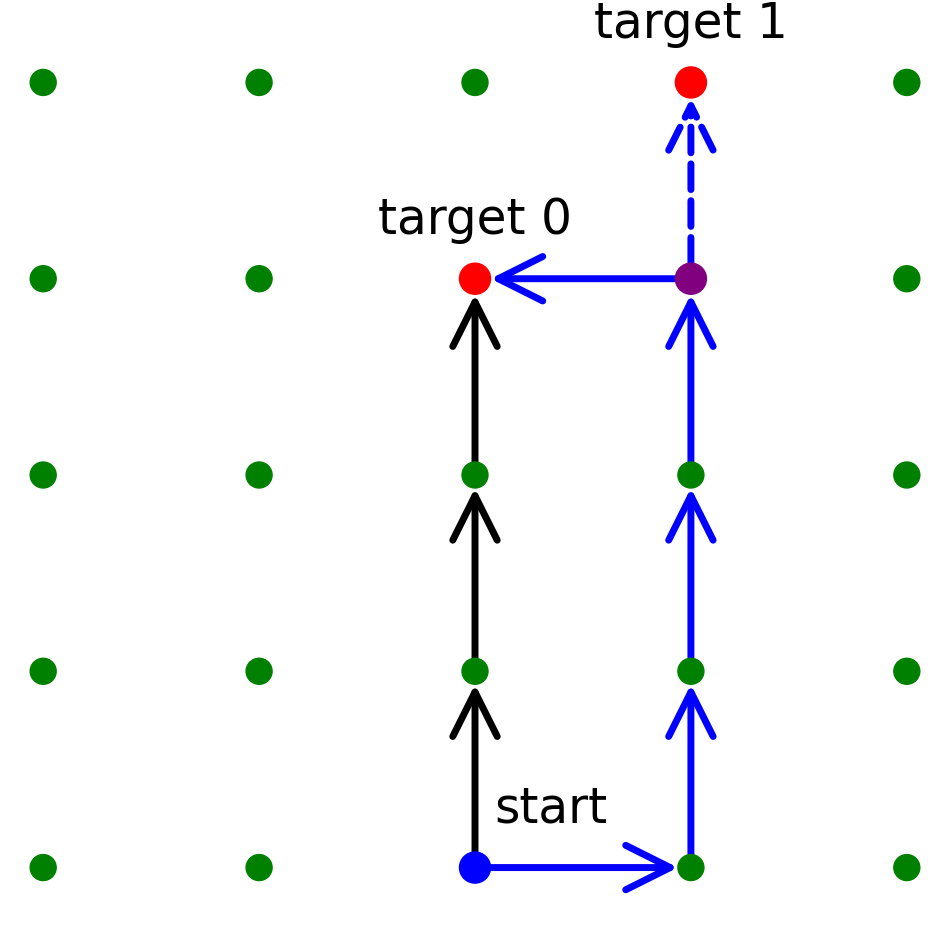}
  \caption{Shortest path (black) versus a longer but deceptive route (blue) to
target 0 (red square).}
  \label{fig:example_shortest}
\end{figure}

\noindent\textbf{Why shortest paths can mislead.}  
Figure \ref{fig:example_shortest} illustrates why simply assuming an
attacker follows the shortest route is insufficient for defense
planning.  The black path is distance–optimal, yet it instantly betrays
the target; the longer blue path hides the agent’s intention until it
leaves the purple node.

\noindent\textbf{Our contribution.}  
We first introduce a \emph{model-agnostic} metric for the
\textbf{counter\-deceptiveness} of a road network, based on the notion
of a \emph{last deceptive point}
\cite{goal_recognition_in_path,deceptive_path_planning}: the farthest
vertex from the start at which at least two targets remain reachable.
Maximising the minimum distance from each target to its last deceptive
point yields the defender’s reaction time against any attacker
behaviour.  We show that the optimal design problem is NP-hard but can
be restricted, without loss of optimality, to the class of
Steiner trees.

Because exact optimization is intractable, we propose
\textsc{Reattachment}, a polynomial-time heuristic that iteratively
detaches the least informative target branch and reattaches it along a
shortest path that improves the metric.  Pre-computed
multiple-source–shortest-path tables and an $\mathrm{A}^*$ fallback keep
each reattachment fast and typically straight/road-parallel.
The algorithm’s complexity is
$\mathcal{O}\!\left(\ell\,|\tau|^{2}\,n\log n\right)$, where $\ell$ is
the number of outer iterations, $|\tau|$ the target count, and $n$ the
grid size. Experiments on thousands of random grids and a realistic Tonopah airport
layout demonstrate consistent improvements—up to 40× in reaction
distance than the best baseline seed tree—under tight budget constraints, while the full run of
\textsc{Reattachment} completes in seconds.

The remainder of the paper formalises the metric
(\Cref{sec: def}), presents theoretical reductions
(\Cref{sec: init_res}), details \textsc{Reattachment}
(\Cref{sec:algo}), and evaluates its performance
(\Cref{sec: res}).

\section{Problem Statement}\label{sec: def}

The objective considered in this paper is to design a road network, between targets in an environment, which is as counterdeceptive as possible under a budgetary constraint. 
The general form of the problem places targets as points on a 2D plane. 
These points must be connected via roads, optimizing counterdeceptiveness subject to budget constraints. 
The set of targets is given \emph{a priori}; therefore every road
network we build must at least connect the start location to each target
so that all targets remain reachable.
Each edge of the graph carries a weight equal to the construction cost of that road segment, and the \emph{budget} is a hard
upper bound on the \emph{sum} of those edge weights.
As we will discuss in \Cref{thrm: nphard}, an exact solution is infeasible to compute.
Thus, we use a approximation technique: generating a grid graph and rounding targets to the graph's vertices. Following this, we use the words graph and (road) network interchangeably.
Our heuristic approximation algorithm will output a highly counterdeceptive graph.
To do this, we first define a \emph{model-agnostic measure of
counterdeception}---the minimum distance an agent can travel before its
goal becomes uniquely determined by the observer---motivated by
\cite{goal_recognition_in_path,deceptive_path_planning}.

Let $\tau\subset\mathbb{R}^{2}$ be a finite set of \emph{targets} and
let $s\in\mathbb{R}^{2}$ be the \emph{start}.  
We work on a directed, weighted graph
\(G=(V,E,w)\), where \(V\) is the node set containing
\(\tau\cup\{s\}\), \(E\) is the set of directed edges (written
\(u\!\to\!v\)), and \(w:E\to\mathbb{R}_{\ge 0}\) assigns a
non-negative length to each edge. A \emph{path} from $u$ to $v$ is a sequence
$u=v_{0}\!\to v_{1}\!\to\cdots\!\to v_{k}=v$; its length is
$w_{P}(u,v)=\sum_{i=0}^{k-1} w(v_{i}\!\to v_{i+1})$.  
We assume every target is reachable from $s$ and $G$ is connected;
otherwise the problem is trivial for unreachable targets.
For any subgraph $S\subseteq G$ we require the same reachability
property, and we define its \emph{cost} as the total edge length
\[
  w(S)=\sum_{u\to v\in E(S)} w(u\!\to v).
\]
Throughout, “distance’’ will refer to these graph–theoretic lengths.

To measure \emph{when} an observer can unmask the attacker’s goal, we
adopt the notion of a \emph{last deceptive point} introduced in
\cite{goal_recognition_in_path,deceptive_path_planning}.
The deceptive agent has a target it is trying to reach, say $t$, chosen from the given set of targets $\tau$.
The agent simultaneously wants to deceive an observer whose goal is to determine which target the agent is heading towards as soon as possible.
Consider some path $P$ the agent is taking from the start locations $s$ to this target $t$.
We define the notion of the \emph{last deceptive point} as the point along $P$ where the observer is completely sure that the agent's true goal is $t$.
\begin{figure}[h]
  \centering
  \begin{subfigure}[t]{0.22\textwidth}
    \centering
    \includegraphics[width=\linewidth]{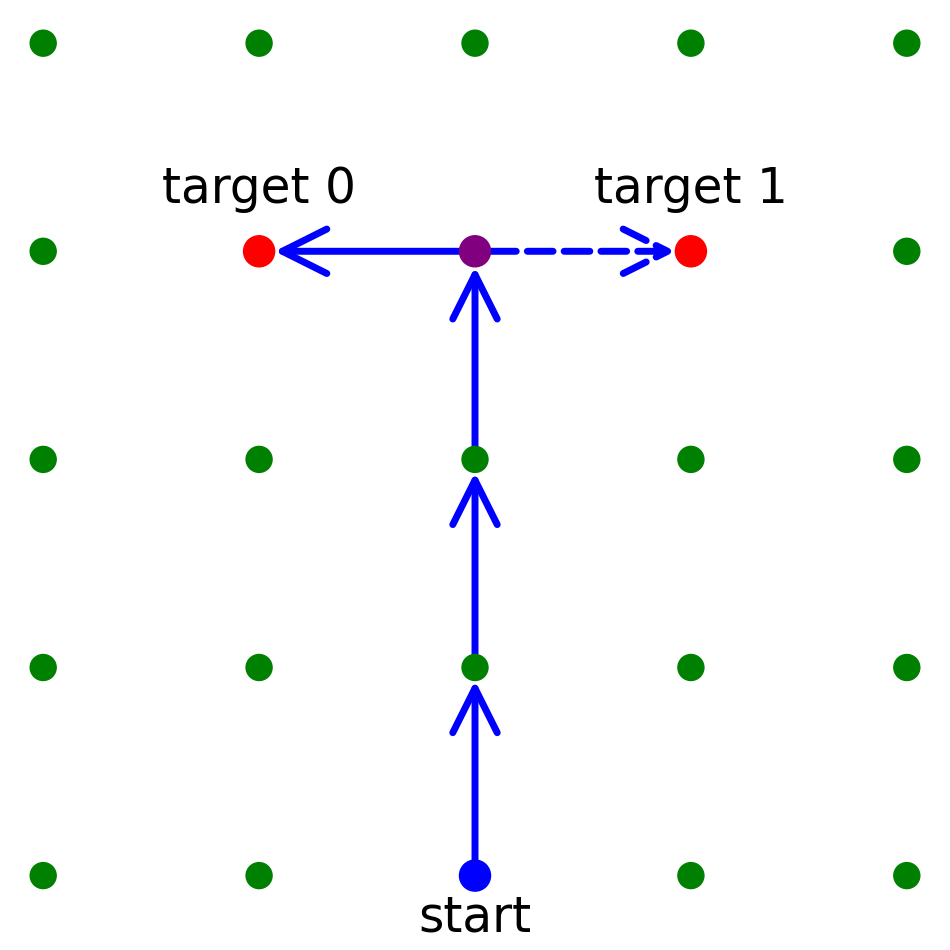}
    \caption{A graph $S$ with a start point (blue) and two targets (red).  
             Leaving the purple node makes the destination certain, so it is
             $l(P,\text{target 0})$ for path $P$ (blue edges).}
    \label{fig: example_last_deceptive}
  \end{subfigure}%
  \hfill
  \begin{subfigure}[t]{0.22\textwidth}
    \centering
    \includegraphics[width=\linewidth]{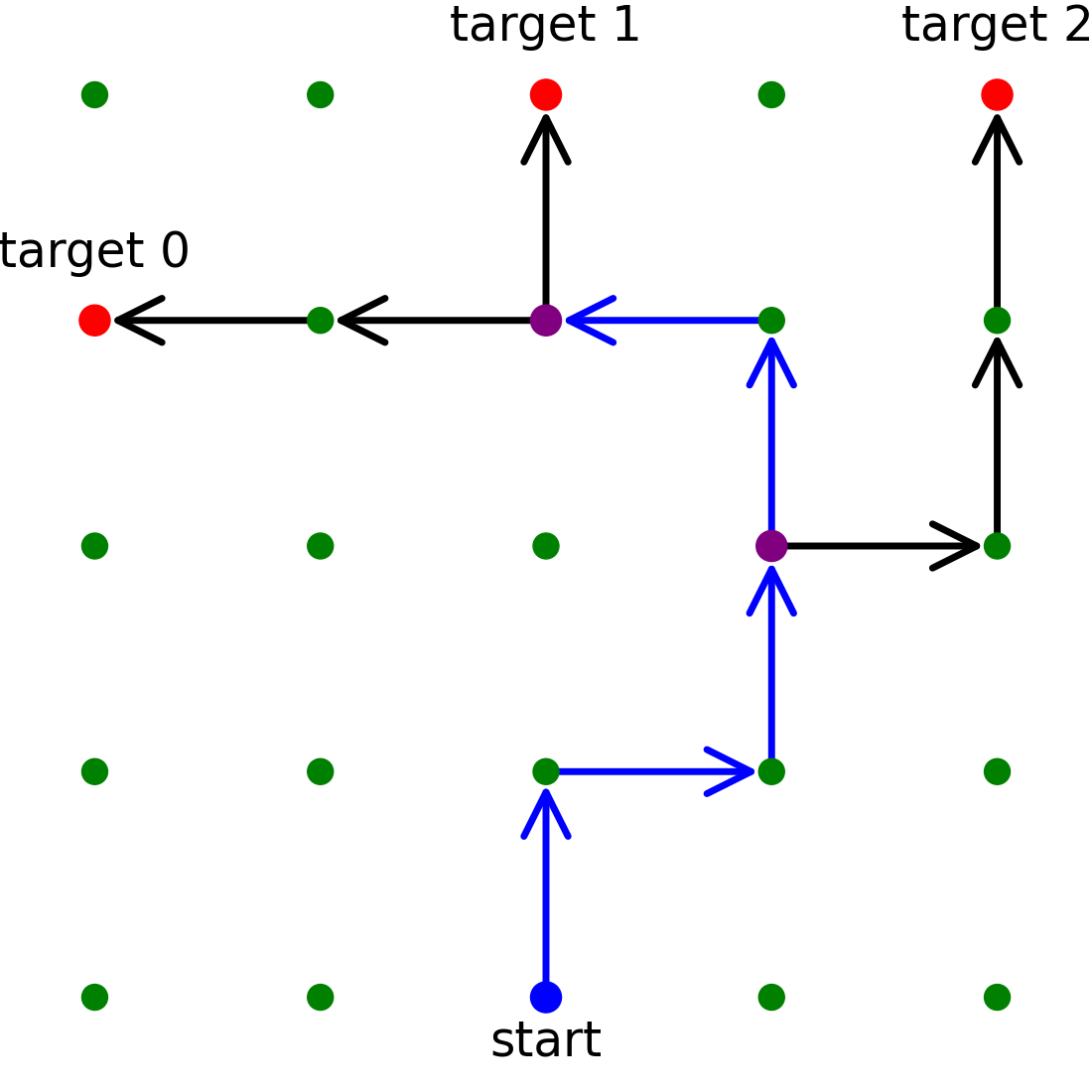}
    \caption{Unique distance illustration. Edges in black are the unique path from $l(P,t)$ to $t$ for each target $t$. Thus, $U_{S}(\text{target 0}) = 2$, $U_{S}(\text{target 1}) = 1$, and $U_{S}(\text{target 2}) = 3$.}
    \label{fig: unique_distance}
  \end{subfigure}
  \caption{Illustrations of last deceptive points (left) and unique
           distances (right).}
  \label{fig:last-and-unique}
\end{figure}
\begin{defn}[Last deceptive point $l(P,t)$]
Let $S$ be a graph with start node $s$ and target set
$\tau$.  Consider a path
$P = v_{0}\!\to v_{1}\!\to\cdots\!\to v_{k}$ in $S$ with
$v_{0}=s$ and $v_{k}=t\in\tau$.
Let $\ell$ be the largest index such that some target
in $\tau\!\setminus\!\{t\}$ is still reachable from $v_{\ell}$ in~$S$.
The vertex $v_{\ell}$ is called the \emph{last deceptive point} of $t$
along~$P$ and is denoted $l(P,t)$.
\end{defn}
Computationally, the reachable-target set at each node $v_i$ on $P$ can be obtained with a single backward traversal, so $l(P,t)$ is easy to locate in practice. \Cref{fig: example_last_deceptive} illustrates the concept.

Intuitively, the \emph{unique distance} of a target is the length of
path an observer has to monitor after the last deceptive point to be
certain of that target.  Formally:
\begin{defn}[Unique Distance, $\textsc{U}_S(t)$]\label{def: uniquedist}
    Let $S$ be a graph with start $s$ and targets $\tau$.
    We define the \textit{unique distance from $s$ to $t \in \tau$}, $\textsc{U}_{S, \tau, s}(t)$, as 
    \begin{align*}
        \textsc{U}_{S, \tau, s}(t)\colon \tau &\to \R_{\geq 0} \\
                                            t &\mapsto \min_{\substack{\text{paths } P \text{ in } S \\ \text{ from } s \text{ to } t}} \Set{w_P\pqty{l(P, t), t}}
    \end{align*} 
    When clear from context, we omit $\tau$ and $s$ and write $\textsc{U}_{S}(t)$.
\end{defn}
See \Cref{fig: unique_distance} for an example of computing $\textsc{U}_S(t)$. For a given target $t \in \tau$, $\textsc{U}_{S}(t)$ is a model-agnostic \emph{model-agnostic} lower bound on how long an observer must wait, in the worst case, before it can be certain the agent is heading for $t$.
Under the reasonable assumption that the observer’s ability to detect or
intercept is uniform across the network, the defender must plan for the worst-case target---one with the \emph{smallest} unique distance.
Because the minimum is taken over \emph{all} \(s\!\to\!t\) paths, this
metric makes no assumptions about the agent’s movement.
With these ingredients we can now formalize the overall
\emph{counterdeceptiveness} of a graph.

\begin{defn}[Counterdeceptiveness, \textsc{CD}]\label{def:metric}
    Given $\tau$, $s$, and $S$ as described before, we define the \emph{counterdeceptiveness} of $S$ as
    \[
        \textsc{CD}_{\tau, s}(S) \defeq \min_{t \in \tau}\Set{\textsc{U}_{S}(t)}.
    \]
    If $\tau$ and $s$ are clear from context we write $\textsc{CD}(S)$.
    For the graph $S$ and the set of targets illustrated in \Cref{fig: unique_distance}, we see that $\textsc{CD}(S) = 1$.
\end{defn}

Returning to the motivational domain of road network design for physical security, if we assume that the observed agent moves at a fixed speed and the observer has perfect observations, \Cref{def:metric} gives a worst case ``reaction time'' available to the observer if the agent is moving towards a different target than the one allowed.
Designing a graph to be as counterdeceptive as possible is now a question of maximizing the value of $\textsc{CD}(S)$ given targets $\tau$ and $s$. 
Inspired by resource constraints in the environment design~\cite{urban_road_network}, we may have a limited budget, and we assume the total resource expenditure is proportional to edge-lengths in the designed graph.
We thus pose the following optimization problem.
\begin{prob}[Budgeted counterdeceptive-design]\label{problem}
Fix a base graph $G=(V,E,w)$, with
start node $s\in V$ and target set $\tau\subset V$.  
Given a budget $b\ge 0$, choose a subgraph
\[
  S\subseteq G\quad\text{s.t.}\quad
  s,\tau\subseteq V(S),\;
  w(S)\le b.
\]
Among all such $S$, maximise the counterdeceptiveness
\(\textsc{CD}_{\tau,s}(S)\):
\[
  S^{\star}\;=\;
  \argmax_{S\subseteq G}
  \Bigl\{\textsc{CD}_{\tau,s}(S)\Bigr\}.
\]
(If a candidate $S$ violates reachability, $\textsc{CD}_{\tau,s}(S)$ is
undefined and the subgraph is deemed infeasible.)
\end{prob}
We will give a heuristic algorithm for solving \Cref{problem}.
However, in order to state such a heuristic, we first give structural optimizations of the types of graphs $S$ that we may consider.

\section{Theoretical Results}\label{sec: init_res}
\begin{figure}[h]
  \centering
  \begin{subfigure}[b]{0.2\textwidth}
    \centering
    \includegraphics[width=\linewidth]{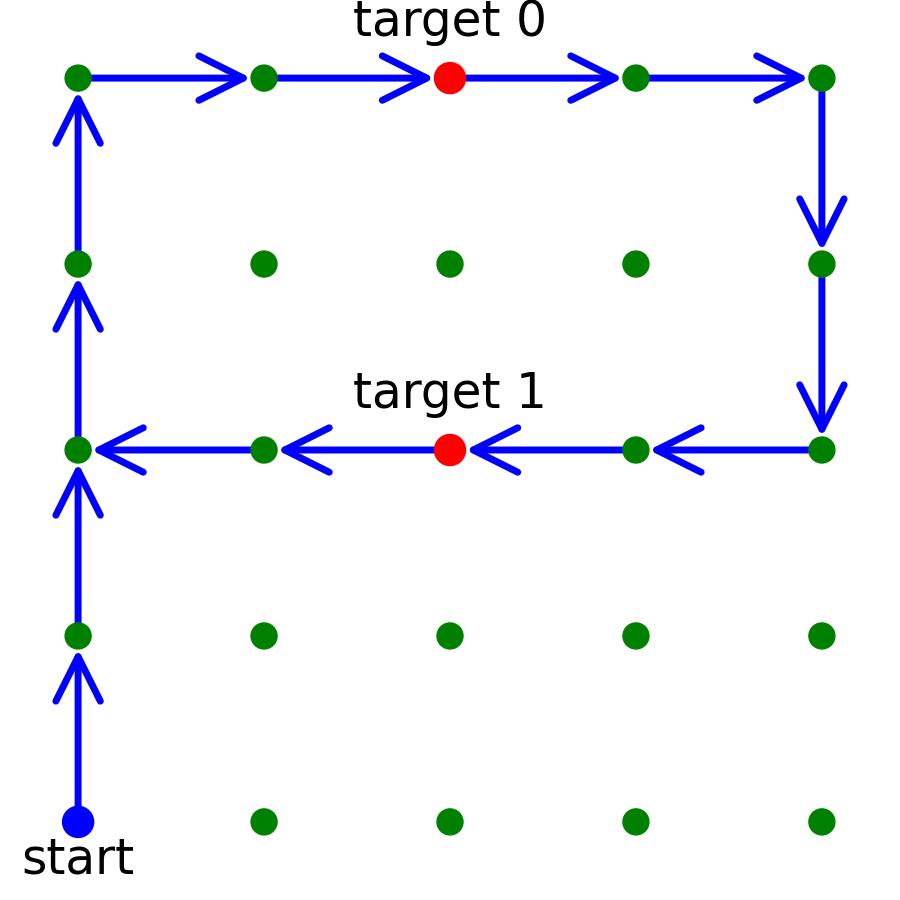}
  \end{subfigure}%
  \hfill
  \begin{subfigure}[b]{0.2\textwidth}
    \centering
    \includegraphics[width=\linewidth]{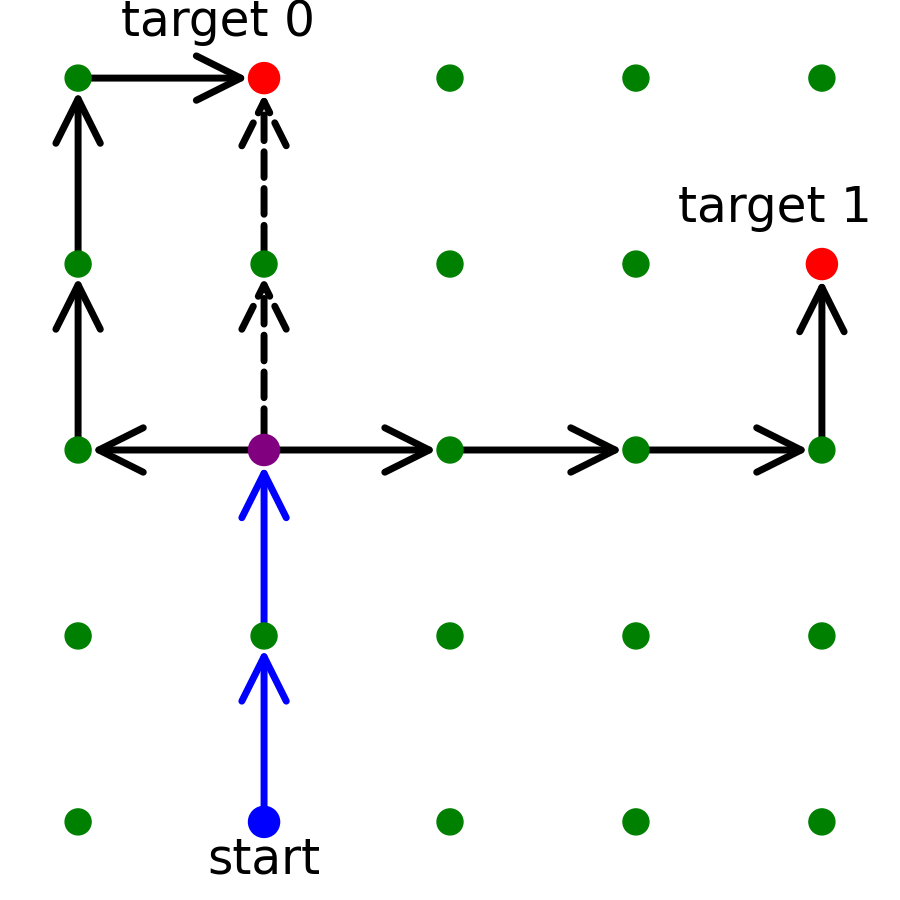}
  \end{subfigure}

  \caption{An example graph with a cycle (left) and an example directed acyclic graph (right). In the right graph, note that the dotted black path can be removed, increasing $\textsc{CD}(S)$ to $4$.}
  \label{fig:cycle}
\end{figure}

Before presenting our heuristic, we first give structural results that shrink
the feasible set of designs in \Cref{problem}.  Recall that we fix a base
graph $G=(V,E,w)$ and seek a subgraph $S\subseteq G$ that spans
$\{s\}\cup\tau$, is connected (each $t\in\tau$ reachable from $s$), and
satisfies $w(S)\le b$.  We show that, without loss of optimality, the
search can be restricted further to \emph{Steiner trees in $G$} spanning
$\{s\}\cup\tau$.
These are trees whose leaves are nodes from a fixed set. In the counterdeception setting, this set will be the targets $\tau$.
Steiner trees are the subject of many combinatorial optimization problems~\cite{min_Steiner, Steiner_survey, Steiner_optimality}; by restricting our search to Steiner trees, we can leverage their well-studied structural properties.
\begin{thrm}\label{thrm: reduction}
Let $\mathcal{S}(\tau,s)$ denote the set of Steiner trees in the base graph $G$
that are rooted at $s$ and have the targets $\tau$ as leaves. Then
\[
  \max_{\substack{S\subseteq G\\ s,\tau\subseteq V(S)\\ w(S)\le b}}
    \textsc{CD}_{\tau,s}(S)
  \;=\;
  \max_{\substack{S\in \mathcal{S}(\tau,s)\\ w(S)\le b}}
    \textsc{CD}_{\tau,s}(S).
\]
In particular, there exists an optimal solution $S^\star\in\mathcal{S}(\tau,s)$.
\end{thrm}
We prove the theorem by transforming any feasible $S$ into a tree of no
smaller counterdeceptiveness and no larger cost.  First, we delete edges
on cycles without reducing $\textsc{CD}$, yielding a directed acyclic
graph.  Next, since agents start at $s$, we prune incoming edges to $s$
and any other sources, leaving a unique source at $s$.  If multiple
$s\!\to\!t$ paths remain for some target $t$, we discard redundant
branches in a way that preserves (or improves) $\textsc{CD}$ and does
not increase weight.  Finally, we trim edges beyond targets so that the
leaves are exactly $\tau$.  Simple examples are shown in \Cref{fig:cycle}. The left graph describes why cycles are nonoptimal. If multiple targets are within a cycle of $S$, then it is impossible to determine which target an adversary is headed towards. The right graph is an example of path discarding to reduce to a tree. These steps show that at least one maximizer is a
Steiner tree; other (non-tree) maximizers may exist only in the case of ties.
\begin{lem}\label{lem: dag}
    At least one maximizer of \Cref{problem} is a directed acyclic graph (DAG).
\end{lem}
\begin{pf}
We will show that given some graph $S$ containing a cycle, we may construct a subgraph $S'$ which is acyclic such that $\textsc{CD}(S') \geq \textsc{CD}(S)$ and $w(S') \leq w(S)$.
By assumption, $S$ is connected and every target is reachable from $s$.
Let $C = c_1 \to c_2 \to \cdots \to c_m \to c_1$ be some cycle in $S$.
If for all targets $t \in \tau$, there was no path from $s$ to $t$ using an edge in $C$ then the lemma is immediate as we may delete any edge from the cycle without decreasing $\textsc{CD}(S)$.
Now suppose that for some target $t \in \tau$ there is a path $P$ from $s$ to $t$ that uses $C$.
After possible relabelling, we may assume that such a path begins with $s \to v_1 \to \cdots \to v_k \to c_1 \to \cdots$.
Form $S'$ by deleting edge $c_m \to c_1$ from $S$.
Since $s$ can reach $c_1$, all nodes on the cycle are still reachable from $s$ meaning that we have not lost the ability to reach any targets from $s$ in $S'$.
Furthermore, deleting an edge cannot create new paths from $s$ to any target, implying that $\textsc{CD}(S') \geq \textsc{CD}(S)$ and $w(S') \leq w(S)$ since all edge weights are non-negative.
Iterating this process forms a subgraph $S' \subseteq S$ with no cycles with $\textsc{CD}(S') \geq \textsc{CD}(S)$ and $w(S') \leq w(S)$.
\end{pf}

\begin{lem}\label{lem: unique_source}
    Define a \emph{source} node in a graph to be a node $u$ where there are no edges $v \to u$.
    Among maximizers of \Cref{problem}, there exists one that is a DAG with the designated start node $s$ as the unique source.
\end{lem}
\begin{pf}
    We will show that given some graph $S$ containing multiple source nodes, we may construct a subgraph $S'$ which has $s$ as the unique source such that $\textsc{CD}(S') \geq \textsc{CD}(S)$ and $w(S') \leq w(S)$.
    We know by \Cref{lem: dag} that any maximizer $S$ of \Cref{problem} is a directed acyclic graph, so without loss of generality let $S$ be such a graph.
    If we have any incoming edges $v \to s$ for some node $v$, then by the acyclic nature of $S$ we have that $v$ is unreachable by any agent starting at $s$.
    Thus, we can remove such edges.
    Similarly, we can remove any sources $u \neq s$ since they are unreachable from $s$ and agents must start at $s$.
    We obtain a graph $S'$ with unique source $s$.
    Since $S$ was acyclic, we have that $S'$ is acyclic and all paths from $s$ to any target $t \in \tau$ are maintained.
    Thus, $\textsc{CD}(S') = \textsc{CD}(S)$ and since edges have been removed, $w(S') \leq w(S)$.
\end{pf}

\begin{lem}\label{lem: tree}
    Among maximizers of \Cref{problem}, there exists one that is a tree rooted at $s$.
\end{lem}
\begin{pf}
    We will show that given some graph $S$ which is not a tree, we may construct a subgraph $S'$ which is a tree such that $\textsc{CD}(S') \geq \textsc{CD}(S)$ and $w(S') \leq w(S)$.
    By \Cref{lem: dag} and \Cref{lem: unique_source} we know that without loss of generality we may take $S$ as a directed acyclic graph with unique source $s$.
    The characterizing difference between directed acyclic graphs and trees is that in directed acyclic graphs, there may be multiple paths between a pair of vertices where in trees there is exactly one path.
    We show that we may remove one of these paths and create a graph with no worse counterdeceptiveness or cost.
    Suppose that a target $t \in \tau$ in $S$ has two different paths $P_1$ and $P_2$ from $s$ to $t$.
    Such paths will share some common vertices at the beginning, then diverge, and then eventually rejoin.
    Without loss of generality, assume that $P_1$ and $P_2$ diverge from each other exactly once and then rejoin exactly once.
    The following argument will work with more notational bookkeeping if this is not the case.
    The paths $P_1$ and $P_2$ take the following form:
    \begin{align*}
        P_1 &= v_0 \to \cdots \to v_i \to v_{x_1} \to \cdots \to v_{x_m} \to v_j \to \cdots \to v_k, \\
        P_2 &= v_0 \to \cdots \to v_i \to v_{y_1} \to \cdots \to v_{y_n} \to v_j \to \cdots \to v_k,
    \end{align*}
    where $v_0 = s$ and $v_k = t$.
    Let $N \defeq \set{v_0, \ldots, v_i, v_{x_1}, \ldots, v_{x_m}, v_{y_1}, \ldots, v_{y_n}, v_j, \ldots, v_k}$ be the set of all nodes in $P_1$ and $P_2$.
    Let $X, Y$ denote the following paths:
    \begin{align*}
        X &= v_{x_1} \to v_{x_2} \to \cdots \to v_{x_m}, \\
        Y &= v_{y_1} \to v_{x_2} \to \cdots \to v_{y_n}.
    \end{align*}
    We consider possible cases of the locations of $l(P_1, t)$ and $l(P_2, t)$.
    \begin{enumerate}
        \item If $l(P_1, t)$ can reach $l(P_2, t)$ (if vice versa, swap $P_1$ and $P_2$), we must have that $l(P_1, t) = l(P_2, t)$.
        This follows from the fact that all targets reachable from $l(P_2, t)$ are also reachable from $l(P_1, t)$.
        To reach $t$ from $s$ we may traverse either $P_1$ or $P_2$.
        If $w_{P_1}(v_{x_1}, v_{x_m}) < w_{P_2}(v_{y_1}, v_{y_n})$ then delete edge $v_{x_m} \to v_j$, otherwise delete edge $v_{y_n} \to v_j$ to form $S'$.
        In both cases, we have disconnected the shorter path, meaning that $\textsc{U}_{S'}(t) \geq \textsc{U}_{S}(t)$.
        Since all nodes in $N$ are still reachable from $s$, we must have that $\textsc{CD}(S') \geq \textsc{CD}(S)$.
        Furthermore, we have that $w(S') \leq w(S)$.
        \item If $l(P_1, t)$ is not reachable from $l(P_2, t)$ and $l(P_2, t)$ is not reachable from $l(P_1, t)$, then we must have that $l(P_1, t) \in X$ and $l(P_2, t) \in Y$.
        Then we may remove either the edge $v_{x_m} \to v_j$ and all targets are still reachable by the agent taking the path through the $v_{y_i}$'s or vice versa.
        Like in the prior case, if $w_{P_1}(v_{x_1}, v_{x_m}) < w_{P_2}(v_{y_1}, v_{y_n})$  then delete edge $v_{x_m} \to v_j$, otherwise delete edge $v_{y_n} \to v_j$ to form $S'$.
        In either case, we have deleted the shorter path, meaning that $\textsc{U}_{S'}(t) \geq \textsc{U}_{S}(t)$.
        Since all nodes are still reachable from $s$, we must have that $\textsc{CD}(S') \geq \textsc{CD}(S)$.
        Furthermore, we have that $w(S') \leq w(S)$.
    \end{enumerate}
    Thus, given two possible paths from $s$ to a target $t$, we have described how to disconnect one path from $S$ to form a new graph $S'$ such that $\textsc{CD}(S') \geq \textsc{CD}(S)$ and $w(S') \leq w(S)$.
    Since $S$ is a finite graph and acyclic, there are only finitely many paths from $s$ to any given target $t$.
    Iterating this process until we are left with unique paths from $s$ to each target $t$ gives a tree $S' \subseteq S$ with $\textsc{CD}(S') \geq \textsc{CD}(S)$ and $w(S') \leq w(S)$.
\end{pf}

We now prove \Cref{thrm: reduction}. \footnote{Complete proofs of
Lemmas~\ref{lem: dag}--\ref{lem: tree} appear in the extended version on
\href{https://arxiv.org/pdf/2311.15093}{arXiv}.}
\begin{pf}
    By Lemmas~\ref{lem: dag}, \ref{lem: unique_source}, and \ref{lem: tree},
the maximum in \Cref{problem} can be achieved by a tree rooted at $s$.
Prune from this tree every edge not lying on an $s\!\to\!t$ path for some
$t\in\tau$; this does not decrease $\textsc{CD}$ and does not increase $w$.
Continue pruning past any target so that the leaves become exactly $\tau$.
Hence an optimal solution exists in $\mathcal{S}(\tau,s)$, and the
maximum over all feasible subgraphs equals the maximum over
$\mathcal{S}(\tau,s)$.
\end{pf}
As a consequence of \Cref{thrm: reduction}, the heuristic algorithm to design an approximate solution to \Cref{problem} will use structural properties of Steiner Trees.
One such property is that in trees, there is a unique path from $s$ to each target $t$. 
Thus, it makes sense to speak of \emph{the} last deceptive point of a target $t$, which we denote as $l(t)$, without referring to a specific path. 

\noindent\textit{Feasibility bound.}
Let $S_{\min}\in\mathcal{S}(\tau,s)$ be a minimum-weight Steiner tree.
If $w(S_{\min})>b$, then \Cref{problem} is infeasible (no subgraph can
connect $s$ to all targets within budget).%
\footnote{Any feasible $S$ must contain a connected subgraph spanning
$\{s\}\cup\tau$, whose cost is at least $w(S_{\min})$.}

Despite the search space reduction proven in \Cref{thrm: reduction}, solving \Cref{problem} is still NP-hard. 
\Cref{thrm: nphard} will briefly show that the complexity of a decision version of \Cref{problem} is \textit{at least} NP-hard.

We can reformulate \Cref{problem} as a decision problem: \emph{Given a set of points $\tau$, minimum counterdeceptiveness $\rho$, and budget $b$, is there a Steiner tree $S \in \mathcal{S}(\tau, s)$ where $\textsc{CD}(S) > \rho$ and $w(S) \leq b$?}. 
In other words, we are only testing for the existence of a tree $S$ that exceeds a counterdeceptiveness of $\rho$ and is under budget $b$. 
The problem we will reduce from is the Euclidean Steiner Minimal Tree (ESMT) problem, which can be described similarly: \textit{Given a set of points $\tau$ and budget $b$, is there a Steiner tree $S$ where $w(S) \leq b$?}

The ESMT problem is known to be NP-hard, but it is unknown if it is NP-complete~\cite{planar_stein}.
We can prove that the decision version of \Cref{problem} is NP-hard by reducing an arbitrary instance of ESMT to an instance of the counterdeceptive problem. 
The consequence of this reduction is that if we had some polynomial time algorithm for finding counterdeceptive road networks, the reduction would give a polynomial time algorithm to solve ESMT.\@
Thus, \Cref{thrm: nphard} proves that the counterdeceptive problem is NP-hard.

\begin{thrm}\label{thrm: nphard}
    Given a set of points $\tau$, minimum counterdeceptiveness $\rho$, and budget $b$, determining if there exists a Steiner tree $S \in \mathcal{S}(\tau, s)$ with $\textsc{CD}(S) > \rho$ and $w(S) < b$ is NP-hard.
\end{thrm}
\begin{pf}
    Take an instance of the ESMT problem with terminal set
$\tau$ and budget $b$: decide whether there exists a Steiner tree of total
weight at most $b$ spanning $\tau$.
    We generate an instance of \Cref{problem}: $s',\tau', \rho', b'$ where $s' = t \in \tau, \tau' = \tau \setminus \{s'\}$, $\rho' = 0$ and $b' = b$. 
    There is a Steiner tree $S$ for $\tau, b$ where $w(S) \leq b$ if and only if there is a Steiner tree $S'$ for $s', \tau', \rho', b'$ where $w(S') \leq b'$. 
    If a tree $S$ for ESMT exists, then $w(S) \leq b'$, and since every target is as a leaf on $S$, $CD(S) > 0$. $S$ is then a valid description for a tree $S'$ optimizing \Cref{problem}, since $\tau = \tau' \cup \{s'\}$. Likewise, if $S'$ exists, then $w(S') \leq b' = b$, thus $S'$ is a valid description for $S$. 
    Our reduction takes constant time, thus the counterdeceptiveness problem is NP-hard.
\end{pf}

When the given budget is close or equal to the cost of the minimum Steiner tree, producing and approximate solution to \Cref{problem} becomes difficult since no solution is allowed to be over budget. 
However, if the budget is reasonable, then finding approximations of counterdeceptive road networks is possible. 
We now discuss the heuristic algorithm to find and improve counterdeceptive road networks, in an effort to optimize \Cref{problem}.

\section{Algorithm}
\label{sec:algo}
\begin{figure}[h]
    \begin{minipage}[c]{0.23\textwidth}
    \centering
    \includegraphics[width=\linewidth]{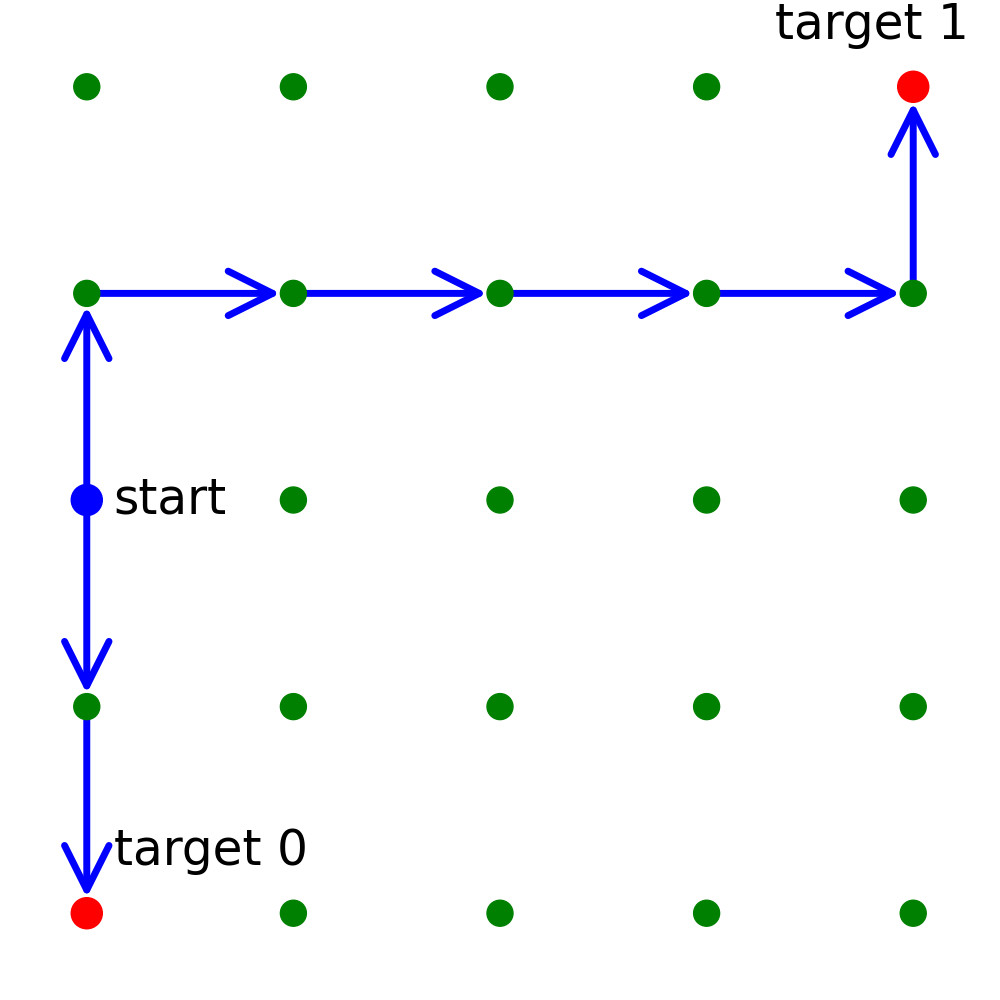}
    \end{minipage}
    \begin{minipage}[c]{0.23\textwidth}
    \centering
    \includegraphics[width=\linewidth]{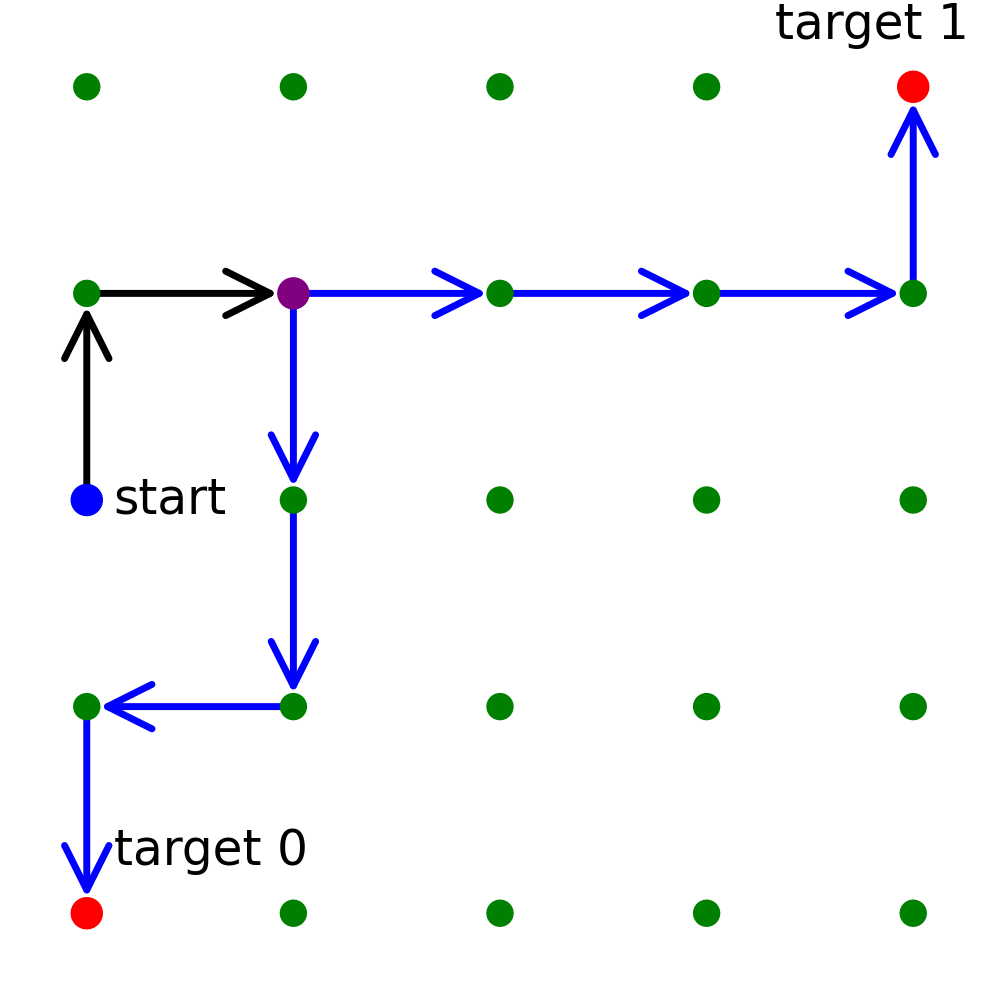}    
    \end{minipage}
    \caption{Before (left) and after (right) of a single reattachment.
    The left graph $S$ has the blue start node as the last deceptive point of each target.
    The right graph $S'$ has the purple point as the last deceptive point of each target.
    Reattachment of the path from start to target $0$ yields an improvement $\textsc{CD}(S') = 4 > \textsc{CD}(S) = 2$.
    }\label{fig: reattach}
\end{figure}

Let targets $\tau$, start location $s$, and budget $b$ be the given inputs to \Cref{problem}.
We introduce a heuristic optimization algorithm, \textsc{Reattachment} (\Cref{alg:reattach}), to search for an approximate solution $S$ to \Cref{problem}.
We know by \Cref{thrm: reduction} that without loss of generality we may take $S$ to be a Steiner tree rooted at $s$ and with the leaves as targets $\tau$, i.e., $S \in \mathcal{S}(\tau, s)$.
The algorithm \textsc{Reattachment} heavily makes use of the fact that in a tree $S$ there is a unique path from $s$ to each target $t \in \tau$.
We can $\textsc{U}_S(t)$ for each $t \in \tau$ and since $\textsc{CD}(S) \defeq \min_{t \in \tau} \set{\textsc{U}_S(t)}$, we identify the least counterdeceptive target $t$ with the lowest value $\textsc{U}_S(t)$.
This leads us to the core idea of the algorithm: \emph{iterative reattachment}. 
Algorithm \textsc{Reattachment} produces better trees by identifying the targets with the lowest $\textsc{U}_S(t)$ and iteratively modifying $S$ to increase this value. 
With some care, we guarantee that improving the counterdeceptiveness of this ``worst'' target does not decrease the overall value of $\textsc{CD}(S)$. 
An example of how reattachment of paths in trees can improve the counterdeceptiveness is shown in \Cref{fig: reattach}.
 
The algorithm \textsc{Reattachment} optimizes a Steiner tree $S \in \mathcal{S}(\tau, s)$ by removing a portion of the tree that connects a single target $t \in \tau$, and finding a new path from $t$ to the remaining portion of the tree, generating $S'$. 
We reject $S'$ if $w(S') > b$, in other words, if the total cost of $S'$ exceeds the given budget.
If $w(S') < b$, we compute $\textsc{CD}(S)$ and compare it to $\textsc{CD}(S')$, accepting $S'$ as the current best tree if $\textsc{CD}(S')>\textsc{CD}(S)$. 
By repeating this procedure until no further improvements, we generate $S_R \in \mathcal{S}(\tau, s)$ that is locally optimal in $\mathcal{S}(\tau, s)$, in the sense that no reattachment of a path in $S_R$ will yield further improvement to $\textsc{CD}(S)$ that is within budget.

We first describe how we initialize \textsc{Reattachment} in \Cref{sec: seed}, where we choose a method to generate a underlying grid graph and starting seed tree to begin iterating from. 
Then, we describe how to compute the counterdeceptiveness of trees in \Cref{sec: counter_indiv}. 
In \Cref{sec: priority}, we discuss how to prioritize targets for reattachment.
We show how to generate new trees via reattachment methods in \Cref{sec: improve_tar}, and describe how to compare the new trees in \Cref{sec: compare_trees}. 
Finally, the whole of \textsc{Reattachment} is detailed in \Cref{sec: final_desc_algo}. 
We will discuss the runtime complexity of each step and show that \textsc{Reattachment} is a polynomial time approximation algorithm.

\subsection{Initialization}\label{sec: seed}

The first step of the approximation algorithm \textsc{Reattachment} is to round the targets and starting point to the nearest nodes on a finite graph, rather than arbitrary points in 2D space.
This allows us to perform reattachment operations quickly, giving us an approximation of the plane.
We choose to use grid graphs of two types as the underlying graph which we will choose edges from.
The use of grid graphs is common in the control and optimization domain~\cite{goal_rec_design_non_obs, goal_rec_design_non_opt, decept_optimal, prob_plan_recognition}.
Furthermore, for the right choice of grids, distances are not distorted by more than a factor of 2 and such grids can be found quickly~\cite{k_enc}.

The first type of grid graph we use are $n_1 \times n_2$ \emph{rectangular grid graphs} with edge weights as distances between the nodes in 2D space.
The second are $n_1 \times n_2$ \emph{triangulated grid graphs}; we take $n_1 \times n_2$ grid graphs and add a central node within each cell with diagonal crossings.
The edge costs are again distances between the nodes in 2D space.
The addition of these diagonal paths creates a reasonable degree of freedom of possible routings for paths while also not adding too many edges.
Since the normal grid graph is a subgraph of the triangulated grid graph, we can also say that distances are not distorted by more than a factor of 2.
Throughout, call the underlying graph $n_1 \times n_2$ grid graph $G$ and take all Steiner trees $S \in \mathcal{S}(\tau, s)$ as subgraphs $S \subseteq G$.
For the purposes of runtime analysis, we will call $n$ the number of nodes in $G$, and since $G$ is planar, the number of edges in $G$ is $O(n)$.

To begin the reattachment procedure, \textsc{Reattachment} also requires some initial tree $S_O \subseteq G$.
A good starting tree is one where it is easy to modify paths from $s$ to the targets, places all the targets on the leaves from the beginning, and is well under budget.
The minimum weight Steiner tree would be meet these criteria, but is infeasible to use as computing the minimal Steiner tree in the setting of planar graphs is NP-Complete~\cite{planar_stein}.
Thus, we opt to use simpler methods in order to generate a seed tree $S_0 \in \mathcal{S}(\tau, s)$.
There are two natural seed trees that we consider: trimmed minimum spanning Steiner trees, and trimmed random spanning Steiner trees.

The first method is to generate a \emph{minimum spanning Steiner seed tree} $S_\textsc{mst}$ by taking a minimum spanning tree of $G$ rooted at $s$ and trimming off any edges that occur past a target that are not connecting another target.
Note that this is not necessarily the minimum weight Steiner tree, as the minimum spanning tree connecting all nodes does not necessarily contain the minimum spanning tree of a subset of nodes.
If $w(S_\textsc{mst}) > b$, we cannot use this method.
An alternative method is to generate a \emph{randomly generated Steiner seed tree} $S_\textsc{rand}$ by generating a random spanning tree of $G$ rooted at $s$ using~\cite{random_span} and performing a similar trimming procedure.
If $w(S_\textsc{rand}) > b$, we may try again, but it is possible that the budget $b$ is too restrictive to generate a valid starting tree. 
Recall that the minimum amount of budget that would generate a feasible solution would be the cost of the minimum Steiner tree. 
Thus, if $b$ is close to such cost, it is highly unlikely that this method would find a valid initial tree $S_O$, and our algorithm will fail.
Note that neither of these methods guarantee that all targets are leaf nodes, only that all leaf nodes are targets.

\subsection{Measuring the Counterdeceptiveness of Individual Targets}\label{sec: counter_indiv}

For each node $u$ in the current tree $S$, we store a value $u.\textit{targets}$ to track how many targets in $\tau$ are reachable from $u$. 
This is a necessary value for computing $\textsc{U}_S(t)$ for each $t \in \tau$ and $\textsc{CD}(S)$ as it allows us to quickly find last deceptive points $l(t)$ for each $t \in \tau$.
The function \textsc{ComputeSuccessors}($S, \tau, s$) computes this value for each node in $S$ in $\bigO{n}$ time using a postorder traversal. 

Our next goal is to compute $\textsc{U}_S(t)$ for each $t \in \tau$. We compute this value directly as defined in \Cref{def: uniquedist}.
\textsc{ComputeSuccessors}($S, \tau, s$) allows us to quickly find $l(t)$ on $S$.
For any node $u$ where $u.\textit{targets}> 1$, but for child $v$ of $u$, $v.\textit{targets} = 1$, let $P_v$ be the path that contains $u\to v$, and terminates with a target $t$.
Then, $l(P_v,t)=u$, as $u$ acts a junction for each child $v$ and the corresponding target.
However, an initial seed tree may place some targets on the tree, such that an adversary is forced to pass through a prior target to reach it. 
Equivalently, these are targets $t \in \tau$ such that $l(t) = t$.
By definition of $\textsc{U}_S(t)$, if $l(t) = t$, then $\textsc{U}_S(t) = 0$. 
We will define such targets as \emph{forced} targets. 
Thus, trees $S$ containing forced targets will have $\textsc{CD}(S) = 0$. 

\begin{figure}[h]
    \centering
    \includegraphics[width=0.25\textwidth]{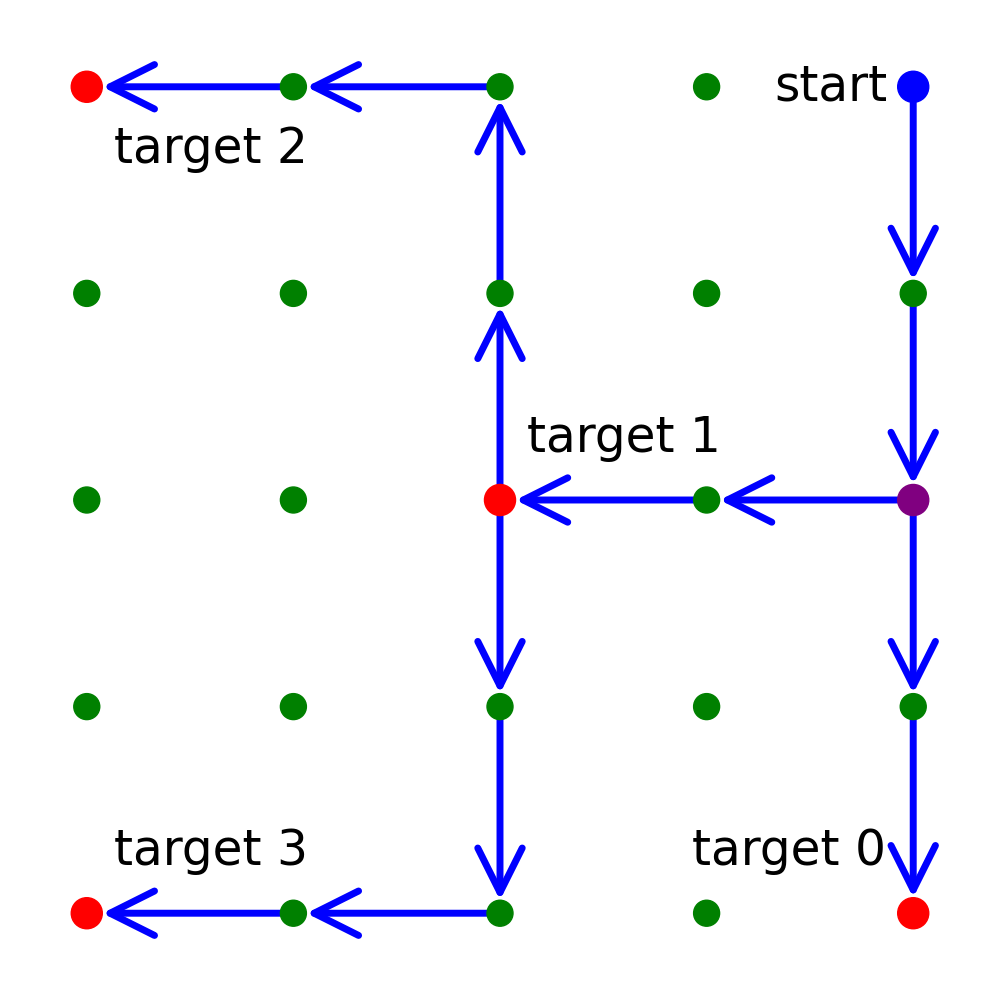}
    \caption{
    Target 1 is forced. 
    To reattach it, target 1 must be completely disconnected from the tree, however, this disconnects target 2 and target 3.
    Thus, reattaching target 1 may not create a graph in $\mathcal{S}(\tau, s)$.
    }\label{fig: forced_reattach}
\end{figure}

The heuristic in \textsc{Reattachment} is to assign a priority to each target $t \in \tau$ such that the least counterdeceptive targets are of higher priority to reattach. 
At every step in \textsc{Reattachment}, we maintain an invariant that our current tree is always a tree $S \in \mathcal{S}(\tau, s)$.
Furthermore, we reattach one path at a time rather than considering multiple reattachments at one time as considering multiple reattachments will make our algorithm slower.
As all forced targets have a counterdeceptiveness of 0, these targets should be among the first to be reattached. 
The only way to reattach a forced target $t$ is to completely remove it from the tree, and make a brand new path. 
However, as shown in \Cref{fig: forced_reattach}, this would remove any paths from $t$ to other targets $t'$, and may disconnect the tree into multiple components, creating a graph $S' \notin \mathcal{S}(\tau, s)$. 
Thus, we cannot naively use $\textsc{U}_S(t)$ as the heuristic to prioritize certain targets to be reattached.
In practice, we assign a custom priority $P_S$ to each target during computation that we use for this purpose.

\subsection{Computing the Priority Heuristic}\label{sec: priority}

\textsc{Reattachment} performs iterative reattachment of nodes in some order.
A decision must be made as to which targets to reattach first, meaning we must assign a \emph{priority} to each target in $\tau$.
If we have forced targets within $S$, then $\textsc{CD}(S) = 0$.
Thus, the main priority is turning each forced target into a leaf node in $S$. 
We can do this by first reattaching the targets that are successors of any forced target in $S$. 
To model this behavior, when a target $t_r$ is reachable from a forced target $t_f$, we set the priority of $t_r$ as $-\textsc{U}_S(t_r)$, otherwise $\textsc{U}_S(t_r)$.
We define a priority $\textsc{P}_S$ as follows:
\begin{defn}[Priority, $\textsc{P}_S(t)$]
    Consider some Steiner tree $S \in \mathcal{S}(\tau, s)$.
    Let the \emph{priority} of a target $t \in \tau$ be defined by a function $\textsc{P}_S$ as follows:
    \begin{align*}
        \textsc{P}_S\colon \tau &\to \R \\
                     t &\mapsto \begin{cases}
                                     \textsc{U}_S(t) & \text{if } l(t) \neq t, \\
                                    -\textsc{U}_S(t) & \text{if } l(t) = t.
                                \end{cases}
    \end{align*}    
\end{defn}

This definition for priority means that for any target $t_r$ reachable from a forced target $t_f$, $\textsc{P}_S(t_r) < \textsc{P}_S(t_f)$ and thus $t_r$ will be reattached before $t_f$. 
\textsc{Reattachment} attempts to reattach targets onto valid locations of the remaining tree. 
Heuristically, we want to remove the paths on the tree that are downstream of some forced target. 
This allows our algorithm to make simpler paths that do not wind around long, unusable branches.
This is why the priority of targets reachable from forced targets is defined to be $-\textsc{U}_S$ since the magnitude of $\textsc{U}_S(t)$ grows with respect to total path length $w_P(l(t), t)$.

Let \textsc{ComputePriority} be the algorithm that computes $\textsc{P}_S(t)$ for each $t \in \tau$ as a sorted list in increasing order. 
Targets with priority $0$ are excluded from the list as these are the forced targets.
As stated in~\Cref{sec: counter_indiv}, detaching a forced target may disconnect the tree, which is not handled by our algorithm.
This list is guaranteed to have some targets of non-zero priority, allowing us to perform at least one reattachment. 
If every target was forced, every target would be able to reach another target in $S$, implying the existence of a cycle.

\textsc{ComputePriority} does a traversal of the tree, taking $\bigO{n}$ time.  
We know from \Cref{sec: counter_indiv} that we can compute in $\bigO{n}$ time the number of reachable targets from reach node in $S$.
From each target, we compute the path from each target to its unique last deceptive point. 
By uniqueness of paths in trees, from each target $t \in \tau$, \textsc{ComputePriority} visits each node at most once.
Thus, $\textsc{ComputePriority}$ runs in $\bigO{n}$ time.

\subsection{Improving Individual Targets}\label{sec: improve_tar}

Now that we can compute the priority $\textsc{P}_S$ of targets to reattach in a tree $S \in \mathcal{S}(\tau, s)$, we can reattach targets in this priority order.
The \emph{invariant} we maintain is that if $S \in \mathcal{S}(\tau, s)$ and we reattach a target $t \in \tau$ forming a new graph $S'$, then $S' \in \mathcal{S}(\tau, s)$.
As stated, we will do this by reattaching one target at a time, never considering multiple reattachments at a time.
Furthermore, we will form $S'$ such that $\textsc{CD}(S') \geq \textsc{CD}(S)$.
We now describe the procedure of reattachment of targets to increase the value of $CD(S)$:

\begin{enumerate}
    \item Let $v \in \tau$ be the target we are reattaching to a new location in $S$.
    \item Mark the branch that connects $v$ to $S$, by following the predecessor nodes of $v$ until we reach a node $x$ such that $x.\textit{targets} > 1$. 
    \item Delete all traversed nodes and edges, except for $v$ and $x$, from $S$ (but not $G$). 
    \item Route the target $v$ back onto a candidate node $c \notin \tau$ via some path in $G$ from $v$ to $c$, creating $S'$.
    We attempt to find a shortest path from $v$ to $c$ using Dijkstra's or $\textsc{A}^*$ as a back up as necessary.
    \item Repeat the previous step trying all possible $c$, keeping the one maximizing $\textsc{CD}(S')$.
\end{enumerate}

To facilitate routing $v$ to $c$, we apply a two stage heuristic.
The first stage is an attempt to connect $c$ to $v$ via a shortest path route though $G$.
We use a Multiple Source Shortest Path (\textsc{MSSP}) algorithm to accomplish this.
Specifically, Dijkstra's algorithm~\cite{dijkstra} is used but various planar-specific MSSP algorithms could be substituted in-place.
These shortest paths are computed as a preprocessing step, and we can recover a path in $\bigO{n}$ time. 
However, if the shortest path between $v$ and $c$ crosses through any point in $S$ and $S'$ would no longer be a tree.
Thus, we use a second stage of reattachment.

In this second stage, we apply a version of the $\textsc{A}^*$ search algorithm~\cite{a_star} to perform path finding around remaining nodes in $S$ directly.
Our path-finding heuristic is the Euclidean distance from $v$ to $c$. 
This execution is akin to running Dijkstra's algorithm, and in the worst case will take $\bigO{n \log n}$ to complete.
Overall, detachment and reattachment of a target node $v$ can be computed in $\bigO{n}$ time and $\bigO{n \log n}$ time by the routines \textsc{DetachTarget} and \textsc{ReattachTarget} respectively.

\subsection{Comparing Trees}\label{sec: compare_trees}

The idea of \textsc{Reattachment} is to iteratively improve trees in order to optimize $\textsc{CD}(S)$ via reattachment procedures described in \Cref{sec: improve_tar}. 
If $\textsc{CD}(S') > \textsc{CD}(S)$, then we naively take the improved $S'$, stopping when no reattachment improves the tree further. 
However, if there are forced targets, or small local minima, relying solely on $\textsc{CD}$ might lead \textsc{Reattachment} to halt before finding meaningful results.
Thus, we introduce secondary criteria which improves factors about the tree not directly represented by $\textsc{CD}(S')$. 
Let $\#F(S)$ be the number of forced targets in $S$ and $A(S) \defeq \frac{1}{\abs{\tau}} \sum_{t \in \tau} \textsc{P}_S(t)$ be the average of all the priorities in $S$. 
We accept a new tree $S'$ compared to $S$ by the following criteria in decreasing priority, noting that we discord $S'$ if $S'$ is ever strictly worse than $S$ in some criterion, i.e., \ if $S' < S$:
\begin{enumerate}
    \item If $\textsc{CD}(S') > \textsc{CD}(S)$;
    \item If $\#F(S') < \#F(S)$;
    \item If $A(S') > A(S)$;
    \item $w(S') < w(S)$.
\end{enumerate}

Our first criterion, $\textsc{CD}(S') > \textsc{CD}(S)$, states we should keep $S'$ as a better tree than $S$ since our primary goal is to maximize $\textsc{CD}$ over the space $\mathcal{S}(\tau, s)$.
Forced targets are the most detrimental targets to the optimization of $\textsc{CD}(S)$.
If any target is forced, then $\textsc{CD}(S) = 0$.
Thus, the first secondary criterion we consider is the number of forced targets since when $\#F(S) = 0$, then $\textsc{CD}(S)$ can be improved in future iterations.

If we can reattach any target other than the highest priority target we may be able to open paths to improve the target with the highest priority. 
By considering the average $A(S)$ of all priority values, this allows us to see improvement to the overall counterdeceptive routes within the tree. 
It is not an undesired consequence to have non-minimal targets be \textit{more} counterdeceptive than previous iterations, however from the wider scope of the problem, the target with minimum counterdeceptiveness determines the baseline quality of the tree. 
If that minimum counterdeceptiveness is equal between reattachments, then improvement in other targets is our next criterion. 
In practice, considering the average priority instead of other similar metrics works well in minimizing the total number of reattachments we make while still improving $S$.

If the averages are equal, we consider the cost as the final criterion. 
If $w(S') < w(S)$, we keep $S'$ as the improved tree. 
Otherwise, the reattachment does not improve $S$ in any of the four criteria, and we keep $S$ and move on to the next candidate node for reattachment.
With an initial graph, tree, priorities, procedures for reattachment, and a concrete way to compare trees in $\mathcal{S}(\tau, s)$, we may describe the whole algorithm \textsc{Reattachment}.

\subsection{The Complete Reattachment Algorithm}\label{sec: final_desc_algo}

We now have the tools necessary to prioritize nodes for reattachment, compute reattachments, and evaluate the quality of different trees.
With all of the prior subroutines in place, we now describe the proposed heuristic reattachment optimization algorithm, \textsc{Reattachment}, in \Cref{alg:reattach}. 
\begin{algorithm}[h]
  \caption{Reattachment$(G,\tau,s,b)$}
  \label{alg:reattach}
  $\,\mathit{paths}\gets \MSSP(G,\tau)\,$\;
  $\,S\gets \SeedTree(G,\tau,s)\,$\;
  \While{reattachment improves $S$}{
    $\,\mathit{priorities}\gets \ComputePriority(S,\tau,s)\,$\;
    $S_{\text{best}}\gets S$\;
    \For{$(t,c)\in \mathit{priorities}$}{
      $S'\gets \DetachTarget(S,t,s)$\;
      \For{$u\in V(S')\setminus V(S)$}{
        $S''\gets \ReattachTarget(G,S',t,u,\mathit{paths})$\;
        \If{$S''>S_{\text{best}}$}{
          $S_{\text{best}}\gets S''$\;
        }
      }
    }
    \uIf{$S_{\text{best}}>S$}{
      $S\gets S_{\text{best}}$\;
      \textbf{break}\;
    }
  }
  \Return{$S$}\;
\end{algorithm}

We combine the asymptotic runtimes given in \Cref{sec: counter_indiv}, \Cref{sec: priority}, and  \Cref{sec: improve_tar} to analyze the runtime of \textsc{Reattachment}.
Running \textsc{Reattachment} takes $\bigO{\ell \cdot \pqty{\abs{\tau}^2 \cdot n \cdot \log n}}$, where $\ell$ is the number of iterations of the outer \textbf{while} loop.
The value of $\ell$ depends greatly on the initial seed tree, the number of targets, and the number of nodes in the graph.
These results will be explored in \Cref{sec: res}.

\section{Experimental Results}\label{sec: res}
In this section, we benchmark our Python~3 implementation of \textsc{Reattachment} against alternative methods of solving \Cref{problem}. All experiments were run on a workstation with an 8th Generation Intel Core i7-8700K
CPU and 32\,GB RAM.
The Python~3 implementation of benchmarks can be found on \href{https://github.com/leadcatlab/counterdeception}{GitHub}\footnote{\href{https://github.com/leadcatlab/counterdeception}{https://github.com/leadcatlab/counterdeception}}.
In particular, we compare \textsc{Reattachment} to brute force methods and fast random generation of Steiner trees. Our evaluation has three components. 
(i) \textbf{Diagnostics:} we quantify how the number of accepted reattachments varies with $|\tau|$ and grid size, and how outcomes depend on the initial seed (trimmed MST vs.\ trimmed random). 
(ii) \textbf{Head-to-head comparison on random grids:} \textsc{Reattachment} vs.\ a time-matched random-Steiner sampler under the same budget. 
(iii) \textbf{Head-to-head comparison on a map inspired by Tonopah airport layout:} the same time-matched, common-budget protocol as in (ii).



We consider rectangular grids $G$ as in \Cref{sec:algo} with given
$\tau$, start $s$, and budget $b$. An exhaustive baseline would enumerate
all spanning trees of $G$, trim non‐target branches to form a Steiner
tree $S\subseteq G$ spanning $\{s\}\cup\tau$, and keep those with
$w(S)\le b$. The count of spanning trees on an $n\times n$ grid grows
exponentially (Sequence A007341 in \cite{oeis}): $4\times4$ has $100{,}352$ trees (we sweep
them in $\approx\!78$\,s on our i7-8700K/32\,GB machine), while
$5\times5$ already has $557{,}568{,}000$---implying $>\!7{,}000$ hours at
the same rate---before even computing \textsc{CD} for each candidate. Hence, we restrict brute force to tiny grids to obtain ground-truth optima; for larger grids, brute force is infeasible and we evaluate \textsc{Reattachment} (with time-matched random baselines) instead.

\subsection{Average Number of Reattachments}
\begin{figure}[h]
  \centering
\includegraphics[width=.49\columnwidth]{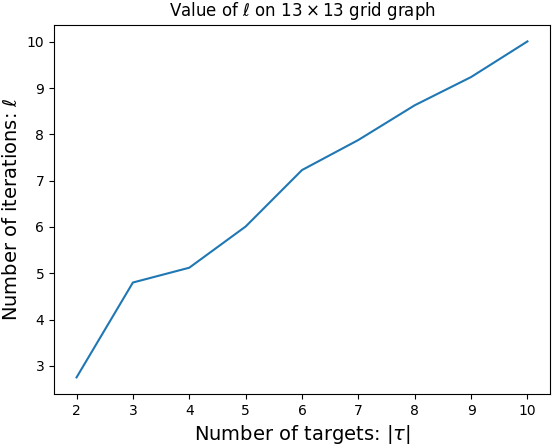}
\includegraphics[width=.49\columnwidth]{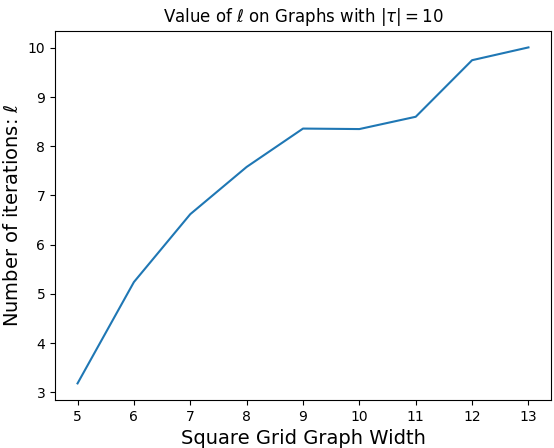}
  \caption{Average value of iterations $\ell$ versus $|\tau|$ on a fixed $13\times13$ grid (left) and $\ell$ versus grid width with $|\tau|=10$ fixed (right).}
  \label{fig:iters}
\end{figure}

Recall that \textsc{Reattachment} runs in $\mathcal{O}\!\left(\,\ell\,|\tau|^{2}\,n\log n\,\right)$ time, where $\ell$ is the number of iterations of the outer while-loop in \Cref{alg:reattach}.
\Cref{fig:iters} describes the average values of $\ell$ on various randomly generated triangulated grid graphs using a randomly generated Steiner spanning tree.
We use an unbounded budget to allow for more possible iterations as limited budgets may rule out possible reattachments causing \textsc{Reattachment} to halt earlier.
In \Cref{fig:iters}, we can see experimentally that $\ell$ grows linearly with respect to both $\tau$ and the width of the grid graphs.
We now compare the quality of the output of \textsc{Reattachment} using different seed trees.
\subsection{Comparison of Different Seed Trees}
We evaluate two initial seeds for \textsc{Reattachment}: a trimmed
minimum–spanning tree (MST) and a trimmed random spanning tree.
On ten $4\times4$ grids (two targets; no budget), we compute the optimal
$\textsc{CD}$ ($\textsc{CD}_{\text{opt}}$) by brute force, then run \textsc{Reattachment} once from
the MST seed and $250$ times from random seeds.  The MST start is
slightly better on these tiny instances as can be seen in Table~\ref{tab:seed-small}, suggesting a $\sim$0.6-approximation.
\begin{table}[h]
  \centering
  \caption{Small grids: mean performance on $10$ random $4\times4$ graphs
  (2 targets, no budget).}
  \label{tab:seed-small}
  \small
  \begin{tabular}{lcc}
    \toprule
    Seed & Runs / graph & Mean $\textsc{CD}/\textsc{CD}_{\text{opt}}$ \\
    \midrule
    Trimmed MST & $1$   & $0.631$ \\
    Trimmed random & $250$ & $0.574$ \\
    \bottomrule
  \end{tabular}
\end{table}

On larger grids, optimal values are unavailable, so we compare seeds by
\emph{win counts}: across all configurations we tally instances where
the MST-seeded run yields higher $\textsc{CD}$ than the random-seeded
run, and vice versa.  Results are mixed (Table~\ref{tab:seed-large}),
indicating no clear winner overall.
\begin{table}[h]
  \centering
  \caption{Larger grids: seed win counts (higher $\textsc{CD}$) by grid size and number of targets.}
  \label{tab:seed-large}
  \small
  \begin{tabular}{cccc}
    \toprule
    Grid & $|\tau|$ & MST wins & Random wins \\
    \midrule
    \multirow{4}{*}{$8\times8$}
      & 2  & 105 & 145 \\
      & 4  & 109 & 141 \\
      & 7  & 143 & 106 \\
      & 10 & 158 &  91 \\
    \midrule
    \multirow{4}{*}{$11\times11$}
      & 2  &  98 & 152 \\
      & 4  & 121 & 129 \\
      & 7  &  95 & 155 \\
      & 10 & 180 &  70 \\
    \midrule
    \multirow{4}{*}{$13\times13$}
      & 2  &  98 & 152 \\
      & 4  & 125 & 125 \\
      & 7  &  88 & 162 \\
      & 10 & 140 & 110 \\
    \midrule
    \textbf{All configs} & — & \textbf{1460} & \textbf{1538} \\
    \bottomrule
  \end{tabular}
\end{table}

\noindent\textbf{Choice going forward.} Since neither seed dominates across settings, we standardize on the
trimmed MST for single–run comparisons (deterministic and slightly
stronger on small grids). When wall-clock permits, we also report a
multi-start variant that samples random seeds. The advantage of considering randomized seed trees will be seen in the next section when considering the framework of running \textsc{Reattachment} multiple times in a fixed time span. 



\subsection{Time-Budgeted Comparison with Random Steiner Sampling}
\label{sec:timebudget}

\textbf{Setup.}
For this benchmark we compare two methods:
(1) \emph{multi-start} \textsc{Reattachment} (Alg.~\ref{alg:reattach}),
run repeatedly with different random seeds and we keep the best
$\textsc{CD}$; and
(2) \emph{random Steiner sampling}, which generates uniformly random
Steiner trees in $G$ and keeps the best $\textsc{CD}$.
For each instance we fix a triangulated grid graph $G$ (as in
\Cref{sec:algo}) with start $s$ and target set $\tau$. Both methods run
on the same $G$, obey the budget constraint, and are scored by
$\textsc{CD}(\cdot)$. We use triangulated grids here to provide richer
routing options (diagonals) for both approaches; earlier brute-force
numbers used smaller rectangular grids only for tractability.

\textbf{Protocol.}
For a fixed wall-clock budget \(T\in\{30,60,90,120,150,180\}\) s, we run:
(i) \emph{multi-start} \textsc{Reattachment} and let \(S_R\) be the best tree found within \(T\); and
(ii) \emph{random Steiner sampling} and let \(S_G\) be the best random
Steiner tree found within \(T\).
Both methods enforce the same budget \(b\).
We repeat on 50 random \(15{\times}15\) triangulated grids with \(|\tau|=8\)
and report win counts (higher \(\textsc{CD}\)).

\begin{table}[h]
\centering
    \begin{tabular}{ccc}
    \toprule
    Time Span & $\textsc{CD}(S_R) > \textsc{CD}(S_G)$ & $\textsc{CD}(S_G) > \textsc{CD}(S_R)$ \\
    \midrule
    30s      & 46                           & 4                            \\
    60s      & 50                           & 0                            \\
    90s      & 50                           & 0                            \\
    120s     & 50                           & 0                            \\
    150s     & 50                           & 0                            \\
    180s     & 50                           & 0                            \\
    \bottomrule
    \end{tabular}
    \caption{Comparison of running \textsc{Reattachment} and randomly generating random spanning trees for fixed time spans.}\label{tab: timespan}
\end{table}

\textbf{Results.}
\Cref{tab: timespan} clearly shows that over a fixed time span, \textsc{Reattachment} almost always produces a more counterdeceptive tree than randomly generating Steiner trees. These results are consistent with the intuition that guided local moves explore useful parts of $\mathcal{S}(\tau,s)$ more efficiently than uninformed sampling.

We now demonstrate the effectiveness of \textsc{Reattachment} on a realistic environment.

\subsection{Design in a Real Facility}

\textbf{Setup.}
The motivating scenario for \textsc{Reattachment} is to design a road network for some commercial or military facility which may consist of many buildings that need to be reachable from some common start which is as counterdeceptive as possible. We use the Tonopah Test Range Airport as a case study.
We block off inaccessible areas (runway, interstate), set the start at the main road off the interstate, and place five targets at buildings.
We construct a triangulated grid (Sec.~\ref{sec:algo}) and remove nodes that fall in blocked regions, yielding a graph with $47{,}293$ nodes and $140{,}384$ edges.

\textbf{Procedure.}
We run \textsc{Reattachment} from a random Steiner seed.
The budget is set to twice the cost of the trimmed MST seed, allowing longer, but potentially more informative, routes while staying realistic.
For a time-matched baseline we also sample as many random Steiner trees as possible.

\begin{figure}[t]
  \centering
  \captionsetup[sub]{font=footnotesize}
  \begin{subfigure}[t]{0.235\textwidth}
    \centering
    \includegraphics[width=\linewidth]{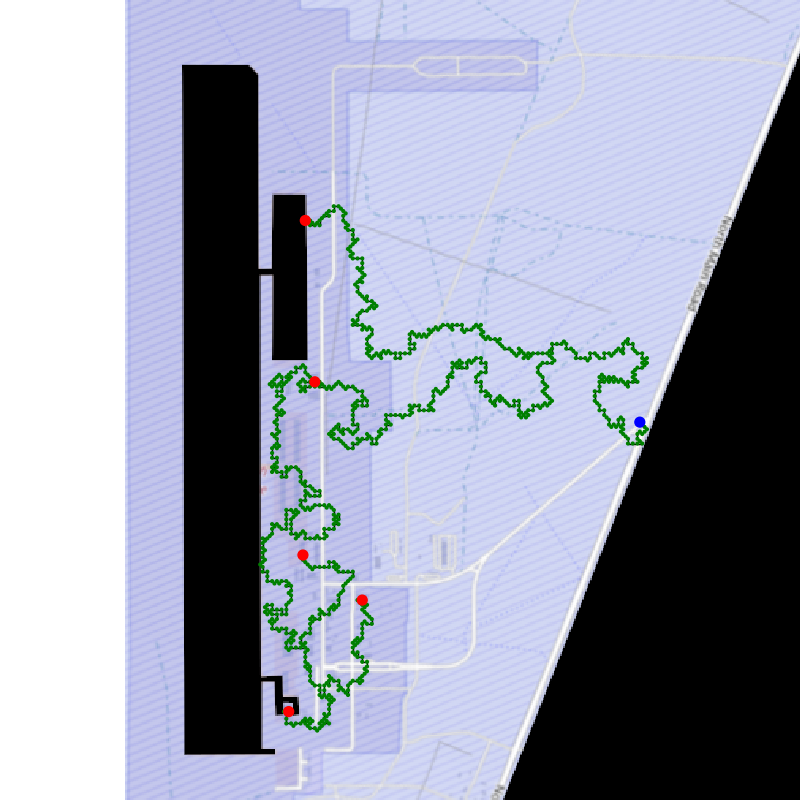}
    \caption{\scriptsize $S_1$: $\textsc{CD}=2.414$, $w=735.131$}
  \end{subfigure}\hfill
  \begin{subfigure}[t]{0.235\textwidth}
    \centering
    \includegraphics[width=\linewidth]{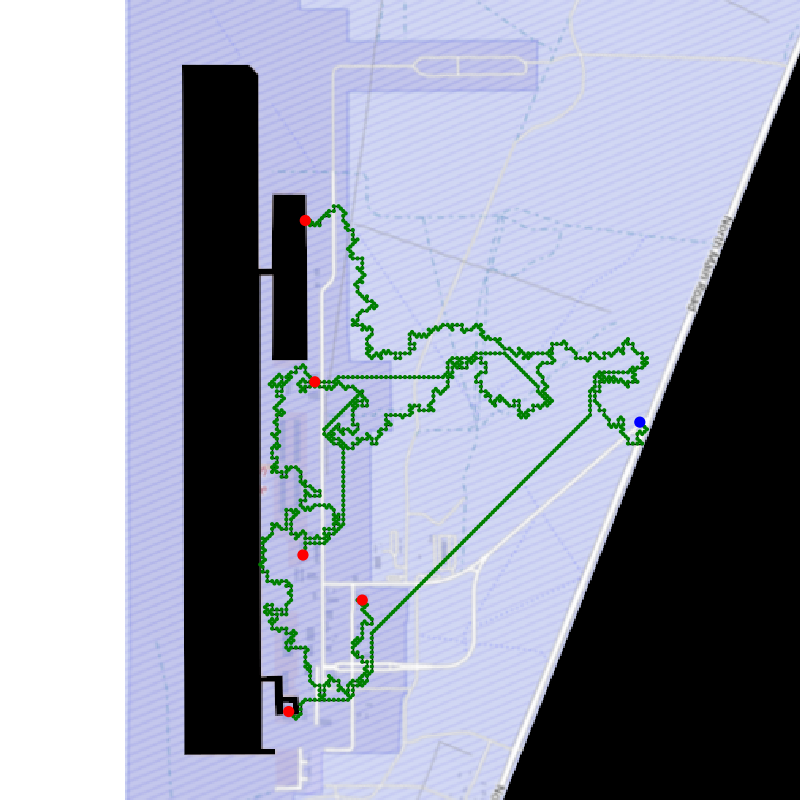}
    \caption{\scriptsize $S_4$: $\textsc{CD}=46.506$, $w=876.670$}
  \end{subfigure}

  \vspace{0.25em}

  \begin{subfigure}[t]{0.235\textwidth}
    \centering
    \includegraphics[width=\linewidth]{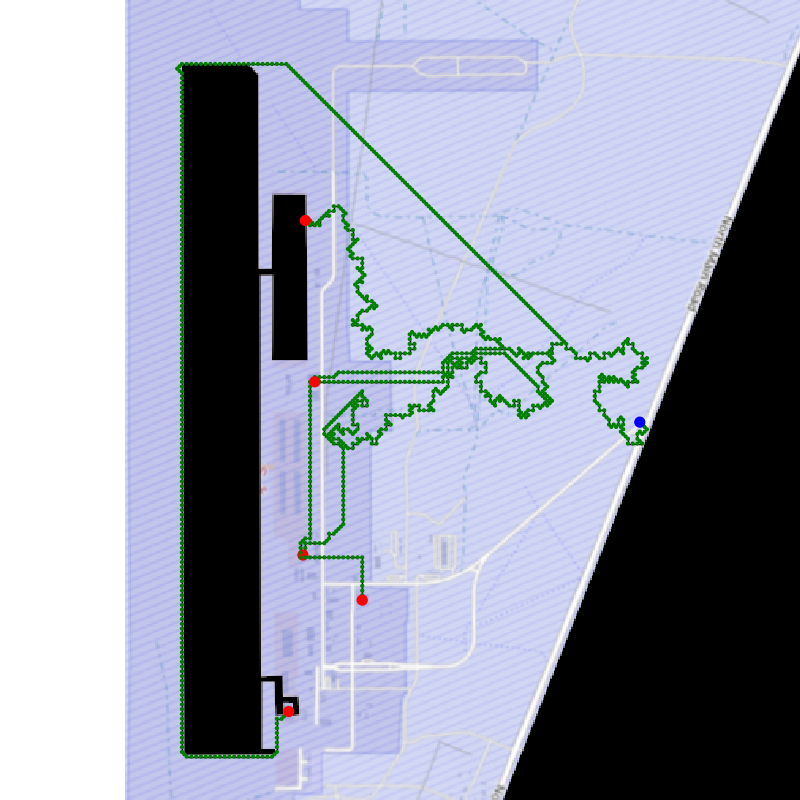}
    \caption{\scriptsize $S_7$: $\textsc{CD}=58.971$, $w=706.734$}
  \end{subfigure}\hfill
  \begin{subfigure}[t]{0.235\textwidth}
    \centering
    \includegraphics[width=\linewidth]{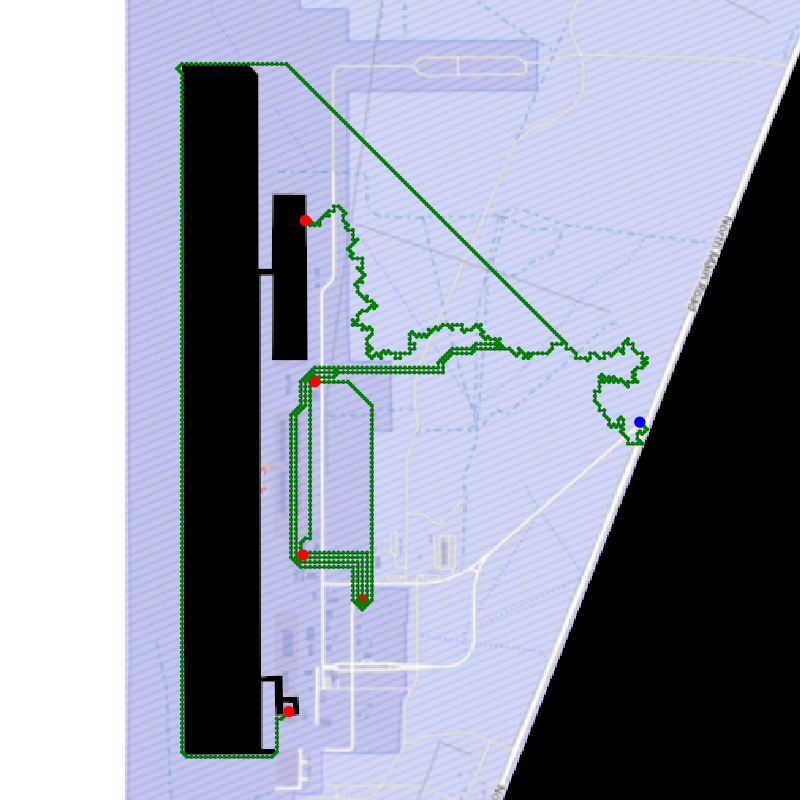}
    \caption{\scriptsize $S_{10}$: $\textsc{CD}=93.314$, $w=876.579$}
  \end{subfigure}

  \caption{A sequence of trees generated by \textsc{Reattachment} starting from an initial random Steiner tree over Tonopah airport with designated start and 5 targets.
    Tree $S_i$ is the $i$-th iteration of \textsc{Reattachment}, and $S_1$ is the initial seed tree while $S_R \defeq S_{10}$ is the final result.}
  \label{fig:airport_gen}
\end{figure}

\textbf{Results.}
\textsc{Reattachment} runs for $1$h $58$m $34$s, performs $10$ reattachments, and outputs $S_R$ with $\textsc{CD}(S_R)=93.314$ and $w(S_R)=876.579$. \Cref{fig:airport_gen} shows representative iterations ($S_1,S_4,S_7,S_{10}$),
with \textsc{CD} rising from $2.414$ to $93.314$ after ten reattachments.
In the same time we generate $4{,}756$ random Steiner trees; the best sampled $S_G$ attains $\textsc{CD}(S_G)=70.406$ with $w(S_G)=891.741$.
Note that $w(S_G)>w(S_R)$ yet $\textsc{CD}(S_G)<\textsc{CD}(S_R)$, i.e., the extra cost in the random baseline does not translate into better counterdeception.

\textbf{Interpretation.}
Randomly sampled trees often contain jagged segments and small kinks that are implausible as roads.
\textsc{Reattachment} tends to insert straighter segments because it routes via shortest paths in the grid.
Our model enforces a tree (unique $s\!\to\!t$ paths), so there are no cross-links between parallel roads; more realistic layout models could relax this with non-tree networks, which is a distinct design problem.

\section{Conclusion and Future Work}
This paper presents a novel model-agnostic metric for network counterdeceptiveness and the associated budget-constrained design problem. The metric uses the last deceptive point---the farthest node from the start at which more than one target remains reachable---and, for each target, measures the \emph{unique distance}, i.e., the path length from that point to the target; maximizing the minimum of these distances yields a guaranteed reaction distance for the defender. The paper proves structural results showing that an optimum can be attained on a Steiner tree (so the search may be restricted accordingly) and establishes NP-hardness of the decision version. To address practical instances, the paper proposes \textsc{Reattachment}, a heuristic that iteratively detaches and reconnects target branches while respecting the budget. Experiments on random grids and a realistic airfield graph demonstrate substantial improvements over baselines within modest runtimes.

\textbf{Applied directions.} Beyond the structural and algorithmic results, several applied directions are immediate.
(i) \emph{Realistic layouts:} add simple geometric regularizers (penalties on sharp turns and very short segments; degree limits; right-of-way masks) so designs resemble plausible roads. 
(ii) \emph{Usability–security trade-off:} formulate a multi-objective variant that maximizes the minimum unique distance while bounding benign users’ travel time and ensuring connectivity. 
(iii) \emph{Cyber analogues:} interpret nodes as services and edges as permitted routes; the “budget’’ would then be rule count, latency, or change cost.

\bibliographystyle{IEEEtran}
\bibliography{IEEEfull}

\end{document}